\documentclass{article}
\usepackage{amsmath, amssymb, latexsym, amscd}
\usepackage{amsfonts}
\author{Christopher A. Francisco\\Department of Mathematics, University of Missouri\\Columbia, MO 65211-4100\\E-mail: chrisf@math.missouri.edu}
\title{Resolutions of small sets of fat points}
\begin{document}
\maketitle

\begin{abstract}
We investigate the minimal graded free resolutions of ideals of at most $n+1$ fat points in general position in $\mathbb P^n$. Our main theorem is that these ideals are componentwise linear. This result yields a number of corollaries, including the multiplicity conjecture of Herzog, Huneke, and Srinivasan in this case. On the computational side, using an iterated mapping cone process, we compute formulas for the graded Betti numbers of ideals associated to two fat points in $\mathbb P^n$, verifying a conjecture of Fatabbi, and at most $n+1$ general double points in $\mathbb P^n$. 
\end{abstract}

\section[1]{Introduction}\label{s:intro}

The study of sets of fat points in projective space is a classical topic that has continued to receive significant attention recently. Many researchers have investigated the Hilbert function, minimal graded free resolution, and other invariants of ideals of fat points, usually with some restriction on the sets of points they consider to make the problems more tractable. In her 2001 paper \cite{Fatabbi}, Fatabbi examines the ideals of at most $n+1$ general fat points in $\mathbb P^n$. The restriction on the number of fat points allows her to work with particularly nice monomial ideals, making it possible to do a thorough analysis of the Hilbert functions, generating sets, and, in some cases, the resolutions, of the ideals of fat points. (See also \cite{tonyi} for methods for finding the Hilbert functions of at most $n+2$ general fat points in $\mathbb{P}^n$.) We extend Fatabbi's work on resolutions in this paper. 

We work in projective space $\mathbb P^n$ and let the corresponding polynomial ring be $R=k[x_0,\dots,x_n]$, where $k$ is a field. Let $(P_0,a_0),\dots,(P_r,a_r)$ be fat points in $\mathbb P^n$, with $a_i$ the multiplicity of the fat point $P_i$. We assume that the $a_i$ are weakly decreasing; that is, $a_0 \ge a_1 \ge \cdots \ge a_r$. Moreover, we assume that $r+1 \le n+1$ so that we have at most $n+1$ fat points in $\mathbb P^n$. Consequently, we may make a change of coordinates, allowing us to suppose that $P_i=[0:\cdots:0:1:0:\cdots:0]$, with the 1 in the $i$-th position, for each $i$. Then the ideal $\mathfrak p_i^{a_i}$ corresponding to the fat point $(P_i,a_i)$ is \[ \mathfrak p_i^{a_i} = (x_0, \dots, x_{i-1},x_{i+1}, \dots, x_n)^{a_i},\] with $x_i$ left out. Thus the ideal of a single fat point is simply a power of an ideal generated by $n$ of the $n+1$ variables. The ideal corresponding to the set of fat points $\{(P_0,a_0),\dots,(P_r,a_r)\}$ is \[ I = \mathfrak p_0^{a_0} \cap \cdots \cap \mathfrak p_r^{a_r}.\] Throughout this paper, $I$ will denote an ideal of this form.

Our paper focuses on the minimal graded free resolutions of the modules $R/I$. In Section 2, we compute the graded Betti numbers of two fat points in $\mathbb P^n$, employing a splitting Fatabbi discovered and used to find the total Betti numbers. The computation proves Fatabbi's conjecture in \cite{Fatabbi} on the graded Betti numbers of two fat points in $\mathbb P^3$. In Section 3, we investigate the free resolution of ideals of at most $n+1$ general fat points with all the same multiplicity (that is, the $a_i$ are all equal). The main result of the paper is Theorem~\ref{mainthm} in Section 4, where we show that if $I$ is the ideal of at most $n+1$ general fat points in $\mathbb P^n$, then $I$ is componentwise linear, a property Herzog and Hibi first introduced in \cite{HerzogHibi}. This has a number of interesting consequences for the graded Betti numbers of $R/I$, which we survey in Section 5. Finally, in the last section, we show that the ideal of $n+2$ fat points in $\mathbb P^n$ may not be componentwise linear, and thus Theorem~\ref{mainthm} is the best we can do.

I dedicate this paper to Graham and Kay Evans on the occasion of Graham's retirement and wish them many more years of happiness.

\section[2]{Two fat points}\label{s:two}

We begin by analyzing the graded Betti numbers of two fat points in $\mathbb P^n$, motivated by Fatabbi's conjecture about the graded Betti numbers of two fat points in $\mathbb P^3$ in \cite{Fatabbi}. We assume throughout that $n \ge 2$ since the $n=1$ case is trivial. The main result in this section, Theorem~\ref{twobetti}, was recently proven independently by Valla \cite{VallaTwo} before we discovered it. We include our proof because our approach is different, and it illustrates the iterated mapping cone technique we shall use throughout the paper. In addition, Fatabbi and Lorenzini have another method for determining the graded Betti numbers of two fat points with different multiplicities (as well as a number of other results on small sets of general fat points in $\mathbb P^n$), and we thank them for sharing their forthcoming paper \cite{FatabbiLorenzini} with us.

By a change of coordinates, we may assume that the ideal corresponding to the two fat points is \[ I = (x_1,\dots,x_n)^{a_0} \cap (x_0,x_2,\dots,x_n)^{a_1},\] where $a_0 \ge a_1$. In \cite{Fatabbi}, Fatabbi finds the minimal generating set of $I$ and then uses a splitting procedure and an inductive argument to compute the total Betti numbers of $I$. She also conjectures formulas for the graded Betti numbers of two fat points in $\mathbb P^3$ and gives some examples as evidence for her conjecture. 

The splitting technique Fatabbi uses first appeared in the paper of Eliahou and Kervaire \cite{EK} in which they describe the minimal free resolution of stable ideals, a class of ideals we shall discuss later. Here we introduce Eliahou and Kervaire's definition of a splittable ideal \cite{EK}.

\newtheorem{definition}{Definition}[section]
\begin{definition}
\label{splitdef}
Let $M$ be a monomial ideal with minimal generating set $G(M)$. $M$ is \textbf{splittable} if there are two nonzero monomial ideals $U$ and $V$ such that:
\begin{enumerate}
\item $G(I)$ is the disjoint union of $G(U)$ and $G(V)$.
\item There is a splitting function from $G(U \cap V)$ to $G(U) \times G(V)$ sending $w$ to $(\phi(w),\psi(w))$ such that:
\begin{enumerate}
\item For all $w \in G(U \cap V)$, $w=\hbox{lcm}(\phi(w),\psi(w))$.
\item If $H$ is a subset of $G(U \cap V)$, then the lcm of $\phi(H)$ and the lcm of $\psi(H)$ both strictly divide the lcm of $H$.
\end{enumerate}
\end{enumerate}
If $U$ and $V$ satisfy the conditions above, we say that $U$ and $V$ are a \textbf{splitting} of $I$.
\end{definition}

The reason this is a useful notion is that we can split an ideal into simpler parts in order to find the minimal free resolution. Fatabbi proves a graded version of a result of Eliahou and Kervaire on the Betti numbers of a splittable ideal \cite{Fatabbi}:

\newtheorem{proposition}[definition]{Proposition}
\begin{proposition}
\label{splitbetti}
\emph{(Eliahou-Kervaire, Fatabbi)} Let $U$ and $V$ be a splitting of a monomial ideal $M$, and let $\beta_{i,j}(R/M)$ denote the $(i,j)$-th Betti number of $R/M$. Then for all $i$ and $j$, \[\beta_{i,j}(R/M) = \beta_{i,j}(R/U) + \beta_{i,j}(R/V) + \beta_{i-1,j}(R/U \cap V).\]
\end{proposition}

Fatabbi shows that the following choices of $U$ and $V$ give a splitting of the ideal $I$ of two fat points: \[ U =(x_2, \dots, x_n)^{a_0}+x_1(x_2,\dots,x_n)^{a_0-1}+\cdots + x_1^{a_0-a_1}(x_2,\dots,x_n)^{a_1} \] \[ V=x_0x_1^{a_0-a_1+1}(x_2,\dots,x_n)^{a_1-1}+x_0^2x_1^{a_0-a_1+2}(x_2,\dots,x_n)^{a_1-2} + \cdots+(x_0^{a_1}x_1^{a_0}) \]

We resolve $R/U$, $R/V$, and $R/U \cap V$ separately and then use Proposition~\ref{splitbetti} to get the graded Betti numbers of $R/I$. Our method uses iterated mapping cone resolutions, and we need a lemma to show that the mapping cone resolutions are minimal at each iteration.

\newtheorem{example}[definition]{Example}
\newtheorem{lemma}[definition]{Lemma}
\begin{lemma}
\label{minmap}
Let $J \subset R=k[x_0,\dots,x_n]$ be a homogeneous ideal with the regularity of $R/J$ at most $d-1$. Let $m$ be a monomial of degree $d$ not in $J$ such that $J:m=(x_{i_1},\dots,x_{i_s})$. Then the mapping cone resolution of $R/(J,m)$ is minimal.
\end{lemma}

\newenvironment{proof}{\noindent {\it Proof}:}{}
\begin{proof}
The lemma follows immediately from Lemma 1.5 in \cite{HerzogTakayama}, using Herzog and Takayama's theory of linear quotients. Alternatively, the result is easy to prove by using the long exact sequence in Tor that the short exact sequence \[ 0 \longrightarrow R/(J:m)(-d) \longrightarrow R/J \longrightarrow R/(J,m) \longrightarrow 0\] induces. \hfill $\square$
\end{proof}

Frequently in the next two sections, we shall need the graded Betti numbers of an ideal of the form $(x_{i_1},\dots,x_{i_s})^d$. We can compute these using Eliahou and Kervaire's resolution of stable ideals \cite{EK} or from the Eagon-Northcott complex.

\begin{lemma}
\label{powerbetti}
Let $M=(x_{i_1},\dots,x_{i_s})^d \subset R=k[x_0,\dots,x_n]$. Then \[ \beta_{q,q+d-1}(R/M) = {d+s-1 \choose d+q-1}{d+q-2 \choose q-1} \] for $q \ge 1$, $\beta_{00}=1$, and all other graded Betti numbers are zero.
\end{lemma}

We begin the process of finding the Betti numbers of $R/I$ by computing the graded Betti numbers of $R/U$.

\begin{lemma}
\label{ubetti} For $q \ge 1$, \[\beta_{q,q+a_0-1}(R/U) = {a_0+n-2 \choose a_0+q-1}{a_0+q-2 \choose q-1} + \sum_{i=1}^{a_0-a_1} {a_0-i+n-2 \choose a_0-i}{n-1 \choose q-1}, \] $\beta_{00}(R/U)=1$, and all other graded Betti numbers are zero.
\end{lemma}

\begin{proof}
We begin with the minimal generators of $(x_2,\dots,x_n)^{a_0}$. The graded Betti numbers of $R/(x_2,\dots,x_n)^{a_0}$ are given by Lemma~\ref{powerbetti}, and it contributes \[ {a_0+n-2 \choose a_0+q-1}{a_0+q-2 \choose q-1} \] to the Betti numbers $\beta_{q,q+a_0}(R/U)$ for $q \ge 1$.

Our goal is to add in one generator of $U$ at a time, computing an iterated mapping cone resolution. The first set of generators we add in with this process are the generators of $x_1(x_2,\dots,x_n)^{a_0-1}$. The order in which we add in generators from $x_1(x_2,\dots,x_n)^{a_0-1}$ will not matter in this case, so we pick descending lex order for specificity. For example, the first ideal we resolve is $(x_2,\dots,x_n)^{a_0}+(x_1x_2^{a_0-1})$. The ideal quotient $(x_2,\dots,x_n)^{a_0}:(x_1x_2^{a_0-1})$ is $(x_2,\dots,x_n)$, which contributes ${n-1 \choose q-1}$ to the graded Betti numbers of $R/U$ since Lemma~\ref{minmap} implies that the mapping cone resolution is minimal. We then add in the generator $x_1x_2^{a_0-2}x_3$ and continue in this way until we have exhausted the generators of $x_1(x_2,\dots,x_n)^{a_0-1}$. Next, we add in the generators of $x_1^2(x_2,\dots,x_n)^{a_0-2}$ in descending lex order and continue this process.

We claim that the ideal quotient at each step is $(x_2,\dots,x_n)$. That is, suppose $U'$ is the ideal \[ U'=(x_2,\dots,x_n)^{a_0} + x_1(x_2,\dots,x_n)^{a_0-1} +\cdots + x_1^{t-1}(x_2,\dots,x_n)^{a_0-t+1} \] \[ + \hbox{ an initial lex segment of generators of } x_1^t(x_2,\dots,x_n)^{a_0-t}.\] Suppose the next monomial to be added into the ideal is $m=x_1^tx_2^{b_2}\cdots x_n^{b_n}$, where the sum of the $b_i$ is $a_0-t$. We need to compute $U':m$. Multiplying $m$ by any of $x_2,\dots,x_n$ lands one in $x_1^{t}(x_2,\dots,x_n)^{a_0-t+1}$, so $(x_2,\dots,x_n) \subseteq U':m$. On the other hand, no multiplication of $m$ by a power of $x_0x_1$ will land in $U'$, and thus $U':m=(x_2,\dots,x_n)$. Again, Lemma~\ref{minmap} implies that the mapping cone resolution is minimal.

Therefore, each minimal generator of $U$ other than the minimal generators of $(x_2,\dots,x_n)^{a_0}$ contributes ${n-1 \choose q-1}$ to the graded Betti numbers of $R/U$. Each $x_1^{i}(x_2,\dots,x_n)^{a_0-i}$ has ${a_0-i+n-2 \choose a_0-i}$ minimal generators; hence the generators other than those in $(x_2,\dots,x_n)^{a_0}$ combine to contribute \[ \sum_{i=1}^{a_0-a_1} {a_0-i+n-2 \choose a_0-i}{n-1 \choose q-1} \] to the graded Betti numbers. \hfill $\square$
\end{proof}

Next, we compute the graded Betti numbers of $V$.

\begin{lemma}
\label{vbetti}
For $q \ge 1$ and $2 \le i \le a_1$, \[ \beta_{q,q+a_0}(R/V)={a_1+n-3 \choose a_1+q-2}{a_1+q-3 \choose q-1}, \] \[ \beta_{q,q+a_0+i-1}(R/V)={a_1-i+n-2 \choose a_1-i}{n-1 \choose q-1}, \] $\beta_{00}(R/V)=1$, and all other graded Betti numbers are zero.
\end{lemma}

\begin{proof}
We resolve $R/V$ the same way we resolved $R/U$. The generators of $V$ in degree $a_0+1$ are the minimal generators of $x_0x_1^{a_0-a_1+1}(x_2,\dots,x_n)^{a_1-1}$. The Betti numbers from this portion of the ideal follow from Lemma~\ref{powerbetti}, and they are the $\beta_{q,q+a_0}(R/V)$ shown above.

To resolve the rest of the ideal, we begin with the generators of $x_0^2x_1^{a_0-a_1+2}(x_2,\dots,x_n)^{a_1-2}$. We can take the generators of this ideal in any order, so we pick descending lex order again. Note that \[ x_0x_1^{a_0-a_1+1}(x_2,\dots,x_n)^{a_1-1}:x_0^2x_1^{a_0-a_1+2}x_2^{a_1-2} = (x_2,\dots,x_n); \] clearly, multiplying by any of $x_2,\dots,x_n$ lands one in $x_0x_1^{a_0-a_1+1}(x_2,\dots,x_n)^{a_1-1}$, but increasing the power of $x_0$ or $x_1$ does not help. By Lemma~\ref{minmap}, the mapping cone resolution of \[ R/(x_0x_1^{a_0-a_1+1}(x_2,\dots,x_n)^{a_1-1}+(x_0^2x_1^{a_0-a_1+2}x_2^{a_1-2}))\] is minimal.

Now let \[ V' = x_0x_1^{a_0-a_1+1}(x_2,\dots,x_n)^{a_1-1}+ \cdots + x_0^{t-1}x_1^{a_0+a_1+t-1}(x_2,\dots,x_n)^{a_1-(t-1)} \] \[ + \hbox{ an initial lex segment of generators of } x_0^tx_1^{a_0-a_1+t}(x_2,\dots,x_n)^{a_1-t}.\] Suppose the next monomial to add into the ideal in the mapping cone process is \[m=x_0^tx_1^{a_0-a_1+t}x_2^{b_2}\cdots x_n^{b_n},\] where the $b_i$ sum to $a_1-t$. We compute $V':m$. It is easy to see that for $i=2,\dots,n$, $x_im \in V'$ because $x_im \in x_0^{t-1}x_1^{a_0+a_1+t-1}(x_2,\dots,x_n)^{a_1-(t-1)}$. Also, $x_0^{b_0}x_1^{b_1}m \not \in V'$ for any choices of $b_0$ and $b_1$, and thus $V':m$ $=(x_2,\dots,x_n)$. 

Therefore the mapping cone resolution is minimal at each iteration by Lemma~\ref{minmap}, and each generator of $V$ in degrees $a_0+2,\dots,a_0+a_1$ contributes ${n-1 \choose q-1}$ to the Betti numbers of $R/V$. Each $x_0^tx_1^{a_0-a_1+i}(x_2,\dots,x_n)^{a_1-i}$ has ${a_1-i+n-2 \choose a_1-i}$ generators, and thus the graded Betti numbers in the statement of the lemma follow. \hfill $\square$
\end{proof}

The last ingredient in computing the graded Betti numbers of $R/I$ is finding the graded Betti numbers of $R/U \cap V$. It is easy to show that \[ U \cap V = x_0x_1^{a_0-a_1+1}(x_2,\dots,x_n)^{a_1}. \] From Lemma~\ref{powerbetti}, we obtain the graded Betti numbers of $R/U \cap V$, which are, for $q \ge 1$, \[ \beta_{q-1,q+a_0}(R/U \cap V)={a_1+n-2 \choose a_1+q-2}{a_1+q-3 \choose q-2},\] $\beta_{00}(R/U \cap V)=1$, and all others zero. We state the graded Betti numbers in this form, using the index $q-1$ for the syzygy, to stay consistent with the formula in Proposition~\ref{splitbetti}.

\newtheorem{theorem}[definition]{Theorem}
\begin{theorem}
\label{twobetti}
Let $I \subset R$ be the ideal of two fat points of multiplicities $a_0 \ge a_1$ in $\mathbb P^n$, where $n \ge 2$. Then for $q \ge 1$, \[ \beta_{q,q+a_0-1}(R/I) = {a_0+n-2 \choose a_0+q-1}{a_0+q-2 \choose q-1} + \sum_{i=1}^{a_0-a_1} {a_0-i+n-2 \choose a_0-i}{n-1 \choose q-1}, \] \[ \beta_{q,q+a_0+i}(R/I)={a_1-i+n-3 \choose a_1-i-1}{n-1 \choose q-1} \hbox{ for } i=0,\dots,a_1-1 ,\] $\beta_{00}(R/I)=1$, and all other graded Betti numbers are zero.
\end{theorem}

\begin{proof}
All the formulas except the one for $\beta_{q,q+a_0}$ follow from Lemmas~\ref{ubetti} and~\ref{vbetti} and Proposition~\ref{splitbetti}. (Note the shift in the range of $i$ from Lemma~\ref{vbetti}.) The formula for $\beta_{q,q+a_0}$ is a consequence of the identity \[ {a_1+n-3 \choose a_1+q-2}{a_1+q-3 \choose q-1} + {a_1+n-2 \choose a_1+q-2}{a_1+q-3 \choose q-2} = {a_1+n-3 \choose a_1-1}{n-1 \choose q-1}.\] \hfill $\square$
\end{proof}

Setting $n=3$ gives the formulas for the graded Betti numbers that Fatabbi conjectured in \cite{Fatabbi}.

\begin{example}\label{twopointex}\emph{
We give the graded Betti numbers of an ideal $I$ of two fat points with multiplicities four and five in $\mathbb P^5$. To display the Betti numbers, we use the notation from the computer algebra system Macaulay 2 \cite{M2}, in which we made all of our computations for this paper. The rows are indexed such that row $d$ contains the Betti numbers $\beta_{i,i+d}$, and $\beta_{i,j}$ is in column $i$ and row $j-i$, where the rows and columns are numbered starting with zero. The graded Betti numbers of $R/I$ are:}

\emph{\begin{tabular}{ccccccccc}
total: & 1 & 126 & 420 & 540 & 315 & 70 &\\
0: & 1 & . & . & . & . & . &\\
1: & . & . & . & . & . & . &\\
2: & . & . & . & . & . & . &\\
3: & . & . & . & . & . & . &\\
4: & . & 91 & 280 & 330 & 175 & 35 &\\
5: & . & 20 & 80 & 120 & 80 & 20 & \\
6: & . & 10 & 40 & 60 & 40 & 10 &\\
7: & . & 4 & 16 & 24 & 16 & 4 &\\
8: & . & 1 & 4 & 6 & 4 & 1 &
\end{tabular}}
\end{example}

Those familiar with the resolutions of stable ideals may notice that this graded Betti diagram looks like it gives the graded Betti numbers of a stable ideal. Recall the definition of a stable ideal: If $m$ is a monomial, let $\max(m)$ be the maximum index of a variable that divides $m$; for example, $\max(x_1^3x_5^3x_7)=7$. We say that a monomial ideal $M$ is \emph{stable} if, whenever $m \in M$ and $i\le \max(m)$, then $\frac{x_i}{x_{\max(m)}}m \in M$. Stable ideals are a generalization of strongly stable ideals, which, when the characteristic of $k$ is zero, are exactly the Borel-fixed ideals. They have received significant attention because generic initial ideals are strongly stable in characteristic zero, stable ideals have convenient combinatorial properties, and Eliahou and Kervaire have computed formulas for the Betti numbers for the graded Betti numbers of stable ideals \cite{EK}. The ideal in Example~\ref{twopointex} is not stable, but its Betti numbers look suspiciously like those of a stable ideal. We shall return to this observation later in the paper.

\section[3]{Fat points with the same multiplicity}\label{s:same}

There are a number of directions one can proceed after analyzing the case of two fat points. One possibility is to investigate sets of more than two but at most $n+1$ general fat points, also requiring that all the fat points have the same multiplicity. Fatabbi discusses an iterative splitting procedure for computing Betti numbers of ideals of at most $n$ general fat points with the same multiplicity in $\mathbb{P}^n$ and examines the first two steps in a splitting in \cite{Fatabbi}. For some recent interesting work on double point schemes using liaison techniques, see \cite{GMS}. We approach the question differently, continuing with the mapping cone idea from Section 2. We show how to use an iterated mapping cone procedure to compute the graded Betti numbers of at most $n+1$ general double points in $\mathbb P^n$. The technique works just as well for general triple points, and we give formulas for those Betti numbers as well. Unfortunately, the bookkeeping gets complicated quickly, and there seems to be no reason to compute explicit formulas for the Betti numbers of ideals of sets of higher order fat points, particularly in light of our main result, Theorem~\ref{mainthm}.

We specialize a result of Fatabbi \cite{Fatabbi} to get the minimal generating set of small sets of general fat points with the same multiplicity.

\begin{proposition}
\label{mingens}
\emph{(Fatabbi)} Let  $I = \mathfrak p_0^{a} \cap \cdots \cap \mathfrak p_r^{a}$ be an ideal of $r+1 \le n+1$ general fat points in $\mathbb P^n$, all with the same multiplicity $a$. Then $I$ is minimally generated by the union of the sets of monomials $G_0,\dots,G_a$, where \[ G_0 = \{ m=x_{r+1}^{b_{r+1}}\cdots x_n^{b_n} \, | \, m \in R_{a}\},\] and for $t=1,\dots,a$, \[ G_t = \{ m=x_0^{b_0} \cdots x_n^{b_n} \in R_{a+t} \, | \, b_i \le t \, \forall \, i=0,\dots,r \hbox{ and } \, \exists \, 0 \le u < v \le r \hbox{ with } b_u=b_v=t\}. \]
\end{proposition}

Thus the degree $a$ generators are all the monomials of degree $a$ involving only the variables $x_{r+1},\dots,x_n$. The higher degree generators in degrees $a+t$ have power of at least two of $x_0,\dots,x_r$ equal to $t$, and no power of $x_0,\dots,x_r$ may exceed $t$. In the case $a=2$, where we have $r+1$ double points, this means that we have minimal generators in degrees 2, 3, and 4, and they have the following form: \[ G_0 = \hbox{degree two monomials in } (x_{r+1},\dots,x_n)^2 \] \[ G_1 = \{x_ix_jx_l \, | \, 0 \le i < j \le r, \, i \not = l \not = j \} \] \[ G_2 = \{ x_0^2x_1^2,\dots,x_0^2x_r^2,x_1^2x_2^2,\dots,x_1^2x_r^2,\dots,x_{r-1}^2x_r^2 \} \] We use this characterization of the minimal generators to compute the graded Betti numbers of at most $n+1$ general double points in $\mathbb P^n$.

\begin{proposition}
\label{doublepoints}
Let $I \subset R$ be the ideal of $P_0,\dots,P_r$, a set of at most $n+1$ general double points in $\mathbb P^n$. Then, for $q \ge 1$, the graded Betti numbers of $R/I$ are \[ \beta_{q,q+1}(R/I)={n-r+1 \choose q+1}{q \choose q-1},\] \[ \beta_{q,q+2}(R/I)= \left (n-r \right ) \sum_{i=0}^{r-1} \left (r-i \right ) {n-i-1 \choose q-1}+ \sum_{i=0}^{r-2} {r-i \choose 2}{n-i-2 \choose q-1},\] \[ \beta_{q,q+3}(R/I)={r+1 \choose 2}{n-1 \choose q-1},\] $\beta_{00}(R/I)=1$, and all other graded Betti numbers are zero.
\end{proposition}

\begin{proof}
The Betti numbers that the monomials in $G_0$ contribute are ${n-r+1 \choose q+1}{q \choose q-1}$ by Lemma~\ref{powerbetti}. To find the contribution of the monomials in $G_1$, we do an iterated mapping cone. This time, we start with the smallest element of $G_1$ in lex order and continue in ascending lex order. (If we pick the largest first, we can get some nonminimal quadratic syzygies that cancel later in the process, and, in particular, we are not in the situation of Lemma~\ref{minmap}.)

We begin by computing $(G_0):(x_{r-1}x_rx_n)$. Multiplication by any of $x_{r+1}, \dots,x_n$ gives a monomial in $(x_{r+1},\dots,x_n)^2$, and we cannot land inside $(G_0)$ by multiplying by any power of $x_0,\dots,x_r$. Thus the ideal quotient is $(x_{r+1},\dots,x_n)$, and because Lemma~\ref{minmap} shows that the mapping cone resolution is minimal, $x_{r-1}x_rx_n$ contributes ${n-r \choose q-1}$ to the graded Betti numbers of $R/I$.

Let $M$ be the ideal generated by $G_0$ and the first $t$ monomials in $G_1$ in ascending lex order. Suppose $x_ix_jx_l$ is the next monomial in $G_1$ to add into the ideal. There are two cases: Without loss of generality, we may assume that either $i$ and $j$ are in $\{0,\dots,r\}$ and $l \in \{r+1,\dots,n\}$, or $i$, $j$, and $l$ are all in $\{0,\dots,r\}$. We wish to compute $M:(x_ix_jx_l)$.

In the first case, where $l \in \{r+1,\dots,n\}$, we may take $0 \le i < j \le r$. Then for all $p \in \{r+1,\dots,n\}$, $x_ix_jx_lx_p \in M$ since $x_lx_p \in G_0$. Thus $(x_{r+1},\dots,x_n) \subset M:(x_ix_jx_l)$. Other variables are also in the ideal quotient, however, because every monomial in $G_1$ less than $x_ix_jx_l$ in lex order is in $M$. Therefore we also have \[(x_{i+1},\dots,x_{j-1},x_{j+1},\dots,x_r) \subset M:(x_ix_jx_l),\] since $x_{i+1}x_jx_l,\dots,x_{j-1}x_jx_l \in M$ and $x_ix_{j+1}x_l,\dots,x_ix_rx_l \in M$. It is easy to see that multiplication by no product of the other variables lands inside $M$. Hence \[ (x_{i+1},\dots,x_{j-1},x_{j+1},\dots,x_n) = M:(x_ix_jx_l).\] By Lemma~\ref{minmap}, the mapping cone resolution of $R/(M,x_ix_jx_l)$ is minimal. There are $(n-j)+$ $(j-i-1)=$ $n-i-1$ variables in the ideal quotient, and by Lemma~\ref{powerbetti}, the ideal quotient contributes ${n-i-1 \choose q-1}$ to the graded Betti numbers of $R/I$.  To compute the Betti numbers arising from the monomials in this case, note that there are $n-r$ choices of $l$, and after $i$ is fixed, there are $r-i$ choices for $j > i$. Since $0 \le i \le r-1$, we have a contribution of \[ \left ( n-r \right ) \sum_{i=0}^{r-1} \left (r-i \right ) {n-i-1 \choose q-1} \] to the graded Betti numbers of $R/I$.

For the other case, we may assume that $0 \le i < j < l \le r$. Again, for all $r+1 \le p \le n$, $x_ix_jx_lx_p \in M$ because $x_ix_jx_p <_{lex} x_ix_jx_l$, and thus $x_ix_jx_p \in M$. Hence $(x_{r+1},\dots,x_n) \subset$ $M:(x_ix_jx_l)$. Next, we consider the variables $x_0,\dots,x_r$, asking which we could multiply by $x_ix_jx_l$ to get a monomial divisible by the elements of $G_1$ in $M$. The elements of $G_1$ in $M$ are all the monomials of $G_1$ less than $x_ix_jx_l$ in lex order, and hence any of $x_{i+1},\dots,x_{j-1},$ $x_{j+1},\dots,x_{l-1},$ $x_{l+1},\dots,x_r$ will multiply $x_ix_jx_l$ into $M$. There is no way to multiply by a product of any of the other variables and land in $G_1 \cap M$, and thus \[(x_{i+1},\dots,x_{j-1},x_{j+1},\dots,x_{l-1},x_{l+1},\dots,x_n) = M:(x_ix_jx_l).\] By Lemma~\ref{minmap} again, the mapping cone resolution of $R/(M,x_ix_jx_l)$ is minimal, and there are $(n-l)+(l-j-1)+(j-i-1)=n-i-2$ variables in the ideal quotient. Lemma~\ref{powerbetti} implies that the monomials in this case add ${n-i-2 \choose q-1}$ to the graded Betti numbers of $R/I$. Once we fix $i$ such that $0 \le i \le r-2$, there are ${r-i \choose 2}$ ways to choose $j$ and $l$ such that $i < j < l \le r$. Hence this case contributes \[ \sum_{i=0}^{r-2} {r-i \choose 2} {n-i-2 \choose q-1} \] to the graded Betti numbers of $R/I$. Combining the two cases gives the graded Betti numbers of the form $\beta_{q,q+2}$.

Finally, we add in the elements of $G_2$. Let $J$ be the ideal generated by $G_0$, $G_1$, and some subset of $G_2$. Pick $0 \le i < j \le r$ such that $x_i^2x_j^2$ is an element of $G_2$ not in $J$. We claim that $J:(x_i^2x_j^2)$ is the ideal generated by all the variables except for $x_i$ and $x_j$. Multiplying $x_i^2x_j^2$ by any $x_l$, where $i \not = l \not = j$, yields a monomial divisible by $x_ix_jx_l$, which is in $G_1$. Additionally, increasing the powers of $x_i$ and $x_j$ on $x_i^2x_j^2$ is no help for getting into $J$. Therefore each element of $G_2$ contributes ${n-1 \choose q-1}$ to the Betti numbers of $R/I$ since the mapping cone resolution is minimal by Lemma~\ref{minmap}. There are ${r+1 \choose 2}$ elements in $G_2$, which gives the formula for $\beta_{q,q+3}(R/I)$. \hfill $\square$
\end{proof}

The same mapping cone technique works for at most $n+1$ general triple points in $\mathbb P^n$. In this case, the ideal has generators in degrees three through six, giving four linear strands in the resolution. We record formulas for the graded Betti numbers to give an idea of what happens going from double to triple points (in particular, the formulas are messier), but we omit the proof, which is much the same as in Proposition~\ref{doublepoints}. It is possible to simplify these formulas somewhat by factoring; we chose to leave them in this form since these expressions are the ones that arise in keeping track of the ideal quotients.

\begin{proposition}
\label{triplepoints}
Let $I \subset R$ be the ideal of $r+1 \le n+1$ general triple points in $\mathbb P^n$. Then, for $q \ge 1$, \[ \beta_{q,q+2}(R/I) = {n-r+2 \choose q+2}{q+1 \choose q-1}, \] \[ \beta_{q,q+3}(R/I) = \sum_{i=0}^{r-1} \left (r-i \right ){n-r+1 \choose 2}{n-i-1 \choose q-1} + \sum_{i=0}^{r-2} {r-i \choose 2} \left (n-r \right ){n-i-2 \choose q-1} \] \[ + \sum_{i=0}^{r-3} {r-i \choose 3}{n-i-3 \choose q-1}, \] \[ \beta_{q,q+4}(R/I)=2{r+1 \choose 3}{n-1 \choose q-1} + {r+1 \choose 3}{n-2 \choose q-1}+(n-r){r+1 \choose 2}{n-1 \choose q-1}, \] \[ \beta_{q,q+5}(R/I)={r+1 \choose 2}{n-1 \choose q-1}. \] All other graded Betti numbers are zero (except $\beta_{00}(R/I)=1$).
\end{proposition}

\begin{example}\label{tripleex}\emph{Consider the ideal $I$ of four triple points in general position in $\mathbb P^5$. The Betti diagram of $R/I$ is:}

\emph{\begin{tabular}{ccccccccc}
total: & 1 & 61 & 203 & 264 & 156 & 35 &\\
0: & 1 & . & . & . & . & . &\\
1: & . & . & . & . & . & . &\\
2: & . & 4 & 3 & . & . & . &\\
3: & . & 27 & 84 & 96 & 48 & 9 &\\
4: & . & 24 & 92 & 132 & 84 & 20 &\\
5: & . & 6 & 24 & 36 & 24 & 6 & 
\end{tabular}}

\emph{Again, the Betti diagram looks like that of a stable ideal, but the ideal is not stable. We devote the next two sections to explaining this phenomenon and its consequences.
}
\end{example}

\section[4]{Fat point ideals and componentwise linearity}\label{s:cwl}

In this section, we prove that ideals of at most $n+1$ general fat points in $\mathbb P^n$ are componentwise linear. This property has many consequences for their resolutions, which we explore in Section 5.

For a homogeneous ideal $J$, let $J_{<d>}$ denote the ideal generated by all the homogeneous elements of degree $d$ in $J$. Herzog and Hibi give the following definition in \cite{HerzogHibi}:

\begin{definition}
\label{cwdef}
Let $J$ be a homogeneous ideal. We call $J$ \textbf{componentwise linear} if $J_{<d>}$ has a $d$-linear resolution for all $d$. That is, for each $d$, $J_{<d>}$ has generators only in degree $d$, first syzygies only in degree $d+1$, etc.
\end{definition}

There are a number of interesting examples of componentwise linear ideals, including stable ideals and Gotzmann ideals (see \cite{HerzogHibi}). Componentwise linear ideals are a natural generalization of ideals with linear resolutions, and their importance first became apparent in a combinatorial application. In \cite{EagonReiner}, Eagon and Reiner prove that a Stanley-Reisner ideal $I_{\Delta}$ associated to a simplicial complex $\Delta$ has a linear resolution if and only if the Alexander dual $\Delta^*$ is Cohen-Macaulay. Herzog and Hibi and Herzog, Reiner, and Welker generalize this result by showing that $I_{\Delta}$ is componentwise linear if and only if $\Delta^*$ is sequentially Cohen-Macaulay, a less restrictive condition than Cohen-Macaulayness that requires a nice filtration of the module $R/I_{\Delta^*}$ in which the quotients are Cohen-Macaulay \cite{HerzogHibi, HRW}.

We begin the process of showing that ideals of at most $n+1$ general fat points in $\mathbb P^n$ are componentwise linear by discussing a notion of Charalambous and Evans \cite{CE:map}.

\begin{definition}
\label{lexwithholesdef}
Let $L$ be a lex ideal in $R=k[x_0,\dots,x_n]$, and let $d_0,\dots,d_n$ be positive integers or infinity. Let $L'$ be the ideal generated by all the minimal generators of $L$ whose degree in $x_i$ is $\le$ $d_i-1$ for all $i$. Then we call $L'$ a \textbf{lex ideal with holes}.
\end{definition}

\begin{example}\emph{
Let $L=(a^3,a^2b,a^2c,ab^3,ab^2c,abc^2) \subset R=k[a,b,c]$. Then $L$ is a lex ideal. Suppose $(d_0,d_1,d_2)=(\infty,3,2)$. We remove all minimal generators of $L$ whose degree in $b$ is 3 or more and whose degree in $c$ is 2 or more. That leaves us with $L'=(a^3,a^2b,a^2c,ab^2c)$, which is a lex ideal with holes.}
\end{example}

We cannot expect the minimal resolution of an arbitrary subideal of a lex ideal to have particularly good properties. However, Charalambous and Evans show that the resolutions of lex ideals with holes do have an especially convenient description \cite{CE:map}.

\begin{theorem}
\label{grahamharathm}
\emph{(Charalambous-Evans)} Let $L$ be a lex ideal, and let $d_0, \dots, d_n$ be positive integers or infinity. Suppose $L'$ is the lex ideal with holes obtained by removing all minimal generators of $L$ whose power of $x_i$ is $\ge d_i$ for each $i$. Then the minimal graded free resolution of $L'$ is a subcomplex of the minimal graded free resolution of $L$. Moreover, one obtains the minimal graded free resolution of $L'$ by deleting all the syzygies in the minimal resolution of $L$ whose degree in $x_i$ exceeds $d_i-1$.
\end{theorem}

We refer the reader to the papers of Eliahou and Kervaire \cite{EK} and Charalambous and Evans \cite{CE:map} for discussions of the basis elements of the syzygy modules (and the degrees of the syzygies) in the minimal free resolution of a lex ideal. We remark only that if a lex ideal with holes does not have too many generators, the process of determining which syzygies survive the deletion process is easy to do by hand.

Let $I = \mathfrak p_0^{a_0} \cap \cdots \cap \mathfrak p_r^{a_r}$ be the ideal of at most $n+1$ general fat points in $\mathbb P^n$ with $a_0 \ge \cdots \ge a_r$ as before. Fatabbi proves in \cite{Fatabbi} that if $t \ge 0$, the set of monomials in $I_{a_0+t}$ is \[ I_{a_0+t} = \{x_0^{b_0}\cdots x_n^{b_n} \in R_{a_0+t} \, | \, b_i \le a_0-a_i+t, i=0,\dots,r\}.\] Thus the powers of $x_0,\dots,x_r$ are restricted, and the powers of $x_{r+1},\dots,x_n$ are not. Using this characterization of the monomials in $I$ in each degree, we show that for each $d$, $I_{<d>}$ is a lex ideal with holes.

\begin{proposition}
\label{arelexholes}
Let $I$ be the ideal of at most $n+1$ general fat points in $\mathbb P^n$. Then for all $t \ge 0$, $I_{<a_0+t>}$ is a lex ideal with holes.
\end{proposition}

\begin{proof}
Let $\mathfrak m=(x_0,\dots,x_n)$. Then $I_{<a_0+t>}$ is generated by all the monomials in $\mathfrak m^{a_0+t}$ except those with power of $x_i$ exceeding $a_0-a_i+t$ for $i=0,\dots,r$. Since $\mathfrak m^{a_0+t}$ is a lex ideal, it follows immediately that $I_{<a_0+t>}$ is a lex ideal with holes. \hfill $\square$
\end{proof}

We can now state our main result.

\begin{theorem}
\label{mainthm}
Let $I$ be the ideal of at most $n+1$ general fat points in $\mathbb P^n$. Then $I$ is componentwise linear.
\end{theorem}

\begin{proof}
We need to show that $I_{<a_0+t>}$ has an $(a_0+t)$-linear resolution for all $t \ge 0$. By Proposition~\ref{arelexholes}, these ideals are lex ideals with holes, so Theorem~\ref{grahamharathm} implies that the minimal resolution of $I_{<a_0+t>}$ is a subcomplex of the minimal resolution of $\mathfrak m^{a_0+t}$, which is $(a_0+t)$-linear. \hfill $\square$
\end{proof}

\section[5]{Consequences of componentwise linearity}\label{s:conseq}

In this section, we discuss the implications of the fat point ideals $I$ being componentwise linear. We begin with a result of Aramova, Herzog, and Hibi \cite{AAH:stable}.

\begin{theorem}
\label{ginbetti}
\emph{(Aramova-Herzog-Hibi)} Let $J$ be a homogeneous ideal in $R=k[x_0,\dots,x_n]$. Let gin$(J)$ be the reverse-lex generic initial ideal of $J$. Then $J$ is componentwise linear if and only if \[ \beta_{i,j}(R/J) = \beta_{i,j}(R/ \hbox{gin}(J)) \] for all $i$ and $j$.
\end{theorem}

Because generic initial ideals are strongly stable in characteristic zero, we get an immediate corollary that explains why the resolutions we examined in Sections 2 and 3 look like those of stable ideals.

\newtheorem{corollary}[definition]{Corollary}
\begin{corollary}
\label{cwgin}
Let $I$ be an ideal of at most $n+1$ general fat points in $\mathbb P^n$, where the underlying field has characteristic zero. Then $I$ has the same graded Betti numbers as a strongly stable ideal, namely its reverse-lex generic initial ideal.
\end{corollary}

We turn now to the question of finding the graded Betti numbers of ideals $I$ of at most $n+1$ general fat points in $\mathbb P^n$. Our goal is to express the graded Betti numbers of $R/I$ in terms of the Betti numbers of the ideals $R/I_{<d>}$; these ideals are lex ideals with holes, and we shall discuss formulas for their Betti numbers later in the section. 

Initially, we note that the graded Betti numbers of componentwise linear ideals satisfy a useful formula of Herzog and Hibi \cite{HerzogHibi}.

\begin{proposition}
\label{cwbetti}
\emph{(Herzog-Hibi)} Let $J$ be a componentwise linear ideal in $R=k[x_0,\dots,x_n]$, and let $\mathfrak m = (x_0,\dots,x_n)$. Then for all $i$ and $d$, \[\beta_{i,i+d}(R/J) = \beta_i(R/J_{<d+1>}) - \beta_i(R/\mathfrak mJ_{<d>}).\]
\end{proposition}

Since the ideals $R/J_{<d>}$ and $R/\mathfrak mJ_{<d>}$ have only linear syzygies, we are writing total Betti numbers for simplicity. The next step is to remove the presence of $\beta_i(R/\mathfrak mJ_{<j>})$ in the formula in Proposition~\ref{cwbetti}, so we determine a formula for its Betti numbers.

\begin{proposition}
\label{mbetti}
Let $J$ be a componentwise linear ideal in $R=k[x_0,\dots,x_n]$. Then for all $i$, \[ \beta_i(R/\mathfrak mJ_{<d>}) = \beta_{i,i+d}(J_{<d>}/\mathfrak m J_{<d>})-\beta_{i+1}(R/J_{<d>}).\]
\end{proposition}

\begin{proof}
We have a short exact sequence \[ 0 \longrightarrow J_{<d>}/\mathfrak m J_{<d>} \longrightarrow R/\mathfrak mJ_{<d>} \longrightarrow R/J_{<d>} \rightarrow 0, \] which induces a long exact sequence of vector spaces in Tor in degree $i+d$: \[ \cdots \rightarrow \hbox{Tor}_{i+1}(k,R/\mathfrak mJ_{<d>})_{i+d} \rightarrow \hbox{Tor}_{i+1}(k,R/J_{<d>})_{i+d} \rightarrow \hbox{Tor}_{i}(k,J_{<d>}/\mathfrak m J_{<d>})_{i+d}\] \[ \rightarrow \hbox{Tor}_{i}(k,R/\mathfrak mJ_{<d>})_{i+d} \rightarrow \hbox{Tor}_{i}(k,R/J_{<d>})_{i+d} \rightarrow \cdots \] The leftmost term is zero since $\mathfrak mJ_{<d>}$ is generated in degree $d+1$. Moreover, the rightmost term is zero because the only nonzero graded Betti numbers of $R/J_{<d>}$, other than $\beta_{00}$, are those of the form $\beta_{i,i+d-1}(R/J_{<d>})$. Thus we have a short exact sequence of vector spaces, and the formula follows. \hfill $\square$
\end{proof}

Finally, we compute the Betti numbers of the modules $J_{<d>}/\mathfrak m J_{<d>}$.

\begin{proposition}
\label{kernelbetti}
Let $J \subset R=k[x_0,\dots,x_n]$ be any homogeneous ideal. Then for $i \ge 0$, \[ \beta_{i,i+d}(J_{<d>}/\frak m J_{<d>})=\beta_{i}(J_{<d>}/\frak m J_{<d>}) = \beta_{1}(R/J_{<d>}) {n+1 \choose i}.\]
\end{proposition}

\begin{proof}
Any degree $d$ element of $J_{<d>}$ multiplied by any $x_i$ lands in $\mathfrak m J_{<d>}$, so we have $n+1$ minimal first syzygies of the form $(0,\dots,x_i,\dots,0)$, $i=0,\dots,n$, for each generator. Any other first syzygy can be written as a combination of these syzygies. The formula for the number of minimal syzygies at each step in the resolution follows immediately. \hfill $\square$
\end{proof}

Combining Propositions~\ref{cwbetti}, \ref{mbetti}, and ~\ref{kernelbetti}, we have formulas for the graded Betti numbers of the fat point ideals in terms of the Betti numbers of the lex ideals with holes $I_{<d>}$.

\begin{theorem}
\label{bettiformulas}
Let $I \subset R$ be the ideal of at most $n+1$ general fat points in $\mathbb P^n$. Then the graded Betti numbers of $R/I$ are given by \[ \beta_{i,i+d}(R/I) = \beta_i(R/I_{<d+1>})+\beta_{i+1}(R/I_{<d>}) - \beta_{1}(R/I_{<d>}) {n+1 \choose i}.\] 
\end{theorem}

In \cite{GHP}, Gasharov, Hibi, and Peeva compute formulas for the graded Betti numbers of lex ideals with holes and, more generally, {\boldmath $a$}-stable ideals. Let $d_i$ be the bounds on the powers of $x_i$ in the lex ideal with holes (so the power of $x_i$ in any minimal generator is less than $d_i$), and for a monomial $m$, define $b(m)=\#\{i \, | \, x_i^{d_i-1}$ divides $m, 0 \le i \le \max(m)-1 \}$. Gasharov, Hibi, and Peeva prove the following theorem.

\begin{theorem}
\label{astable}
\emph{(Gasharov-Hibi-Peeva)} Let $L$ be a lex ideal with holes in $R=k[x_0,\dots,x_n]$ with all its generators in degree $d$. Then the graded Betti numbers of $R/L$ are \[ \beta_{i,i+d-1}(R/L) = \sum_{m \in G(L)} {\max(m)-b(m) \choose i-1} \] for $i \ge 1$, and $\beta_{00}(R/L)=1$, with all other graded Betti numbers zero.
\end{theorem}

We have adjusted the formula from \cite{GHP} to reflect that we are working with variables $x_0,\dots,x_n$ instead of $x_1,\dots,x_n$, and we have restricted the theorem to the case we need. Corollary 2.3 in \cite{GHP} is much more general, and the Eliahou-Kervaire formulas for the Betti numbers of stable ideals actually follow from that result.

As a consequence of Theorem~\ref{astable}, we can, in principle, get formulas for the graded Betti numbers of any ideal of at most $n+1$ general fat points in $\mathbb P^n$. Given the multiplicities of the points, we can use Fatabbi's characterization to list the monomials in $I$ in each degree. Theorem~\ref{astable} allows us to write down the graded Betti numbers of the $I_{<a_o+t>}$ without any difficult computation, and then we can apply Theorem~\ref{bettiformulas} to compute the graded Betti numbers of $R/I$.

Our final application is to a conjecture of Herzog, Huneke, and Srinivasan on the multiplicity of a polynomial ring modulo a homogeneous ideal. 

\newtheorem{conjecture}[definition]{Conjecture}
\begin{conjecture}
\label{multconj}
\emph{(Huneke-Srinivasan, Herzog-Srinivasan)} Let $J$ be a homogeneous ideal of codimension $c$ in $R=k[x_0,\dots,x_n]$ such that $R/J$ is Cohen-Macaulay, and let $e(R/J)$ be the multiplicity of $R/J$. Let $m_i$ be the minimal degree of a syzygy at step $i$ in the minimal graded free resolution of $R/J$, and let $M_i$ be the corresponding maximum. Then \[ \frac{1}{c!} \prod_{i=1}^c m_i \le e(R/J) \le \frac{1}{c!} \prod_{i=1}^c M_i.\]
\end{conjecture}

Conjecture~\ref{multconj} is known in a number of special cases but is open in general, even for monomial ideals in codimension three and above. Recently, there has been interest in proving it for configurations of points in $\mathbb P^n$; see, for example, the paper of Gold, Schenck, and Srinivasan \cite{GSS}. The fact that ideals of at most $n+1$ general fat points in $\mathbb P^n$ are componentwise linear gives the result for free in that case.

\begin{proposition}
\label{multpoints}
Let $I \subset R$ be the ideal of at most $n+1$ general fat points in $\mathbb P^n$ over a field of characteristic zero. Then $R/I$ satisfies Conjecture~\ref{multconj}.
\end{proposition}

\begin{proof}
R{\"o}mer proves Conjecture~\ref{multconj} in \cite{Roemer} for all componentwise linear ideals over a field of characteristic zero, noting that the result follows directly from Theorem~\ref{ginbetti} since Conjecture~\ref{multconj} is true for stable ideals. \hfill $\square$
\end{proof}

We note that the upper bound is not hard to show regardless of characteristic. The upper bound for the ideal $I$ of a set of at most $n+1$ general fat points in $\mathbb P^n$ is \[ \frac{1}{n!} \prod_{i=0}^{n-1} (a_0+a_1+i) \] because there are generators in degree $a_0+a_1$, and $R/I$ is Cohen-Macaulay. Since $\dim R/I=1$, the multiplicity cannot exceed that of $k[x_1,\dots,x_n]/(x_1,\dots,x_n)^{a_0+a_1}$, which is equal to the upper bound.

\section[6]{Larger sets of fat points}\label{s:larger}

A natural question is whether we can extend the results of Section 4 to sets of more than $n+1$ general fat points in $\mathbb P^n$. The following example shows that Theorem~\ref{mainthm} does not hold for $n+2$ fat points in $\mathbb P^n$.
\begin{example}\label{n+2}\emph{Consider four double points in general position in $\mathbb P^2$. We can take the ideal defining these fat points to be \[ I=(b,c)^2 \cap (a,c)^2 \cap (a,b)^2 \cap (a-b,a-c)^2 \] in $R=k[a,b,c]$. The minimal graded free resolution of $R/I$ is \[ 0 \rightarrow R^2(-6) \rightarrow R^3(-4) \rightarrow R \rightarrow R/I \rightarrow 0. \] Clearly, the resolution of $I_{<4>}$ is not 4-linear, and thus $I$ is not componentwise linear.}
\end{example} 

There are a number of directions in which one can proceed from here. First, while Example~\ref{n+2} shows that the ideal of $n+2$ general fat points in $\mathbb P^n$ will not necessarily be componentwise linear, special arrangements of points or fat points may yield componentwise linear ideals. It would be interesting to investigate ideals corresponding to various geometric objects; we are confident there are more ideals arising from geometry that are componentwise linear. Additionally, it would be particularly useful to have more tests for componentwise linearity available to aid in checking for the condition. Finally, the natural long-term goal in this area is the question of Herzog, Reiner, and Welker: Suppose $M$ is a sequentially Cohen-Macaulay module. Does there exist a natural dual module $M^*$ that has a componentwise linear resolution? If so, this would generalize the theorem of Herzog and Hibi on sequentially Cohen-Macaulay simplicial complexes.

\providecommand{\bysame}{\leavevmode\hbox to3em{\hrulefill}\thinspace}

\end{document}